# COXETER GROUPS, 2-COMPLETION, PERIMETER REDUCTION AND SUBGROUP SEPARABILITY

Paul E. Schupp

ABSTRACT. We show that all groups in a very large class of Coxeter groups are locally quasiconvex and have uniform membership problem solvable in quadratic time. If a group in the class satisfies a further hypothesis it is subgroup separable and relevant homomorphisms are also calculable in quadratic time. The algorithm also decides if a finitely generated subgroup has finite index.

## §1. Introduction.

Several years ago the author [S2] raised the question of whether or not small cancellation methods could be used to investigate questions about finitely generated subgroups of "sufficiently nice" groups— in particular the solvability of the membership problem and the Howson property. Rips [R] replied "No" by constructing, for every metric small cancellation condition $C'\left(\frac{1}{n}\right)$, finitely presented groups which satisfy the condition but which have unsolvable membership problem and are neither Howson nor coherent. However, recent work of McCammond and Wise [M-W] and of Arzhantseva and Olshanskii [A-O] shows that the answer is actually "Yes" in many cases. McCammond and Wise introduced the use of "distributive" small cancellation hypotheses where the condition involves how the generators are distributed among all the defining relators and the ingenious idea of "perimeter reduction" where one counts "what is missing".

In this paper we introduce the idea that if one has a "suitable" subgroup graph, $\Delta_1(H)$, of a subgroup H, one can construct from it a "complete" subgroup graph, $\Delta_2(H)$ which has most of the properties of subgroup graphs in the case of free groups and which directly reveals desired information about $H$. . Thus there are two distinct stages. The first is to show that finitely generated subgroups indeed have "suitable" subgroup graphs. In this paper we show that for a very large class of Coxeter groups, and for surface groups, one can start with the usual subgroup graph $\Delta_0(H)$ for an arbitrary finitely generated subgroup $H$ and use perimeter reduction to obtain a subgroup graph $\Delta_1(H)$ which satisfies a certain "Relator Path Property", which is a precise definition of "suitable".

We then show that any subgroup graph whatsoever with the Relator Path Property can be completed to its "2-completion," $\Delta_2(H)$. A "2-complete" subgroup graph has the desired strong properties. Indeed, one can construct the 2-completion of any graph( not necessarily a subgroup graph), having the Relator Path Property and we shall use this fact later. Also, Kapovich and Schupp[K-S] show that the Arzhantseva-Olshanskii technique mentioned above, which is quite different from perimeter reduction, also constructs subgroup graphs with the Relator Path Property.

The graph $\Delta_2(H)$ describes the "Dehn hull" of $H$. That is, any Dehn reduced word $w$ ( one not containing more than half of a defining relator) represents an element of $H$ if and only if $w$ is the label on a closed path at the basepoint of $\Delta_2(H)$. The graph $\Delta_2(H)$ is the transition graph of a finite automaton which accepts a Dehn reduced word $w$ if and only if $w \in H$. Since any geodesic word which is Dehn reduced, this immediately shows that $H$ is rational and quasiconvex and that the diameter of $\Delta_2(H)$ is a quasiconvexity constant for $H$. Next, the subgroup $H$ has finite index if and only if the graph $\Delta_2(H)$ is *full,* that is, there are edges labelled by all generators incident at all vertices.

Perhaps most interestingly, in the case of Coxeter groups, if $G$ satisfies an additional "Separability Condition " then $G$ is subgroup separable and, given $w \notin H$, by slightly modifying $\Delta_2(H)$ to $\Delta_2(H,w)$, one can directly read off from the latter graph a homomorphism $\varphi$ from $G$ into a finite symmetric group such that $\varphi(w) \notin \varphi(H)$. Indeed, there is a uniform quadratic time algorithm calculating all the above information.

A *Coxeter group* $G$ is a group with a finite presentation
$$G = \langle a_1, \ldots, a_n; a_i^2, \quad (a_i a_j)^{m_{ij}}, \ i \neq j \rangle$$
where the $m_{ij} > 1$ and we may have $m_{ij} = \infty$, signifying the absence of a relation between $a_i$ and $a_j$. The *Coxeter matrix* of the presentation is the symmetric matrix $M = (m_{ij})$ where each $m_{ii} = 1$. The *modified Coxeter graph* $\Gamma$ of the presentation is the graph with vertices $a_1, \ldots, a_n$ and, if $m_{ij} < \infty$, there is an edge between $a_i$ and $a_j$ labelled by $m_{ij}$. Note that our modified graph differs from the usual Coxeter graph in that we omit edges labelled by $\infty$ and include edges labelled by 2. Thus the absence of an edge between $a_i$ and $a_j$ means that $m_{ij} = \infty$ When we say that we are "given" a Coxeter group, we always mean that we are given such a finite presentation.

In general, Coxeter groups may be very bad with respect to the properties in which we are interested. Let $B = F_2 \times F_2$ be the direct product of two free groups of rank 2. Mikhailova's Theorem [see L-S] shows that there is a fixed finitely generated subgroup $H$ of $B$ such that the membership problem for $H$ in $B$ is unsolvable. Furthermore, $B$ is neither Howson nor coherent. The class of Coxeter groups is closed under both free products and direct products and the free product $K$ of three cyclic groups of order two contains free subgroups of rank two. Thus $K \times K$ is a six-generator right-angled Coxeter group containing a copy of $B$ and so has unsolvable membership problem and is neither Howson nor coherent.

We turn to formulating a suitable Reduction Hypothesis which will define the class of Coxeter groups to which our results apply. First of all, we assume that there are at least three generators and that all $m_{ij} \geq 4$. In the terminology of Appel and Schupp [A-S], we only consider groups of *extra-large type*. Such groups satisfy the small cancellation hypothesis $C'\left(\frac{1}{7}\right)$ and thus are hyperbolic in the sense of Gromov. (The only use we make of small cancellation theory is in the statement of two facts.)

The set $A$ of generators is ordered

$$a_1 < a_2 < \cdots < a_n$$

by increasing subscript. This induces the shortlex or canonical ordering in the set $A^*$ of all finite words on the generators. Write $|u|$ for the length of $u$ and define $u < v$ if $|u| < |v|$ or if $|u| = |v|$ and $u$ precedes $v$ in the lexicographical ordering on $A^*$ induced by the ordering on $A$. The *normal form* of an element $g$ of $G$ is the least element $u$ in the shortlex ordering such that $u = g$ in $G$.

Tits [T] proved that for any Coxeter group, an arbitrary word can be transformed into normal form by a finite sequence of operations:

(I) Cancel an occurrence of a generator and its inverse (which is just the generator in the case of Coxeter groups).
(II) Replace an occurrence of half a defining relator by the inverse of the other half.

Since we also want to state results about surface groups at the end of this article, we point out that the same result holds for suitable small cancellation groups.

**Normal Form Theorem.** *Let $G = \langle X; R \rangle$ be a finitely presented group in which the symmetrized set $\hat{R}$ generated by $R$ satisfies the small cancellation $C'\left(\frac{1}{7}\right)$, has all pieces of length 1 and has only elements of even length. Then the set $N$ of shortlex normal forms of elements of $G$ is a regular language and there is a linear time algorithm which, when given an arbitrary word $w$ calculates the normal form of $w$. A word $w$ can be transformed into normal form by a finite sequence of operations I and II.*

This theorem is essentially now well-known but we include a proof along with the discussion of Dehn convexity in the Appendix.

We want to study finitely generated subgroups by constructing their subgroup graphs exactly as in the case of a free group. Given a Coxeter group $G = \langle A; R \rangle$, a *graph over* $G$ is a graph in which edges are labelled by generators from $A$. If the subgroup $H$ is generated by $\{h_1, \ldots, h_m\}$ we start with a bouquet of $m$ loops arranged around a basepoint $O_H$. The $i$-th loop is then subdivided into edges labelled by generators so that one reads $h_i$ counterclockwise around the loop. (The $h_i$ are arbitrary words, not necessarily in normal form.) Call this graph $\Delta_0(H)$.

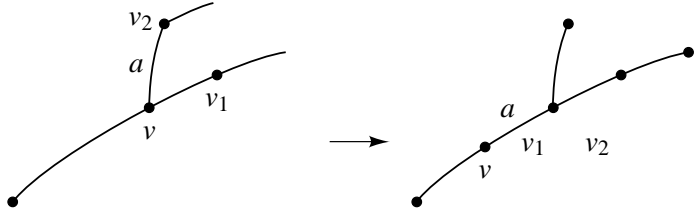

Figure 1

Since the generators have order two, we read the same label in both directions along an edge. If two edges with the same label $a$ are incident at a vertex $v$ we fold (identify) the two edges into a single edge.

Folding is the only operation on the graph which is needed in the case of a free group and takes care of the relations $a_i^2 = 1$ in the case of Coxeter groups, but we now need to deal with the other relations. Two generators $a_i$ and $a_j$ are *related* if $m_{ij} < \infty$, that is, there is a relator $(a_i a_j)^{m_{ij}}$. We view relators as cyclic words. Since the generators have order 2, a relator and its inverse are the same cyclic word. For each $i$, let $\rho_i$ denote the number of other generators $a_j$, $j \neq i$, to which $a_i$ is related. Set $\rho_{ij} = \max\{\rho_i, \rho_j\}$. We shall see later that the following hypothesis is sufficient to use McCammond-Wise perimeter reduction to construct $\Delta_1(H)$.

**The Reduction Hypothesis..** $G = \langle A; R \rangle$ *is of extra-large type and there is a subset $C \subseteq A$ such that every defining relator $(a_i a_j)^{m_{ij}}$ contains a generator from $C$ and satisfies the following condition.*

(1) *If both $a_i$ and $a_j$ are in $C$ then $m_{ij} > \dfrac{3}{2}\rho_{ij}$,*
(2) *If $a_i \in C$ and $a_j \notin C$ then $m_{ij} > 2\rho_i$.*

One can always, of course, make the "uniform choice" $C = A$ and then the hypothesis is that all $m_{ij} > \dfrac{3}{2}\rho_{ij}$.

The Reduction Hypothesis is best viewed as a condition on the modified Coxeter graph $\Gamma$ of $G$. Since the vertices $v_i$ correspond to the generators $a_i$, $\rho_i$ is the degree of the vertex $v_i$. Since $m_{ij} < \infty$ labels an edge $e_{ij}$ between $v_i$ and $v_j$, $\rho_{ij}$ is the maximum of the degrees of the endpoints of $e_{ij}$. A *vertex cover* of the graph $\Gamma$ is a set $C$ of vertices such that every edge has at least one endpoint in $C$. The Hypothesis is that $m_{ij} > \frac{3}{2}\rho_{ij}$ if both endpoints of $e_{ij}$ are in $C$ and $m_{ij} > 2\rho_i$ if only one endpoint $v_i \in C$. For example, suppose that $\Gamma$ is a central triangulation of a polygon, say with $a_n$ the central vertex. In this case, taking $C = \{a_1, \ldots, a_{n-1}\}$ imposes the condition $m_{ij} \geq 5$ on the outer edges and the condition $m_{ni} \geq 7$ on the central edges, and this condition is independent of $n$. In general, the condition on $m_{ij}$ depends only on degrees of the endpoints of the edge $e_{ij}$.

If we used 0 instead of $\infty$ in the Coxeter matrix to represent the absence of a relation, then having a large number of missing relators is a sparseness condition on the matrix. A sparse matrix allows relatively small values for the $m_{ij}$.

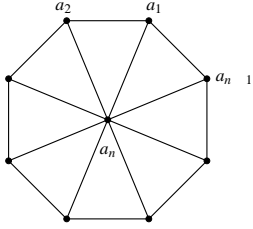

Figure 2

The additional hypothesis which we need to establish the Coxeter group $G$ is subgroup separable is easily stated. Why it is needed will be clear when we discuss the structure of the graphs $\Delta_2(H)$.

**The Separability Condition..** *All $m_{ij}$ are even and if $m_{ij}$ is the label on an edge which forms part of a triangle in the modified Coxeter graph $\Gamma$ then $m_{ij}$ is also divisible by $3$ (and is therefore divisible by $6$).*

Note that if the Coxeter graph $\Gamma$ does not contain any triangles then the hypothesis needed for subgroup separability is the Reduction Hypothesis plus having all $m_{ij}$ even.

We can now state our main theorem.

**Theorem I.** *Then there is a fixed quadratic time algorithm which, when given a Coxeter group $G = \langle A; R \rangle$ satisfying the Reduction Hypothesis, a tuple $(h_1, \ldots, h_m)$ of generators for a subgroup $H$ and an element $w$ of $G$:*

(1) *Calculates the graph $\Delta_2(H)$ which, in particular, is the graph of a finite state automaton which accepts a Dehn reduced word $u$ if and only if $u \in H$.*
(2) *Calculates a quasiconvexity constant for $H$.*
(3) *Decides if $H$ has finite index.*
(4) *Decides if $w \in H$.*
(5) *If $w \notin H$ and $G$ satisfies the Separability Condition, the algorithm explicitly writes out a homomorphism $\varphi$ from $G$ to a finite symmetric group such that $\varphi(w) \notin \varphi(H)$.*

Before proceeding with the proof we make a few comments. (See also Gersten and Short [G-S] who pointed out the close connection between rationality and quasiconvexity.) Since we have fixed shortlex normal forms, we say that a subgroup $H$ is *rational* if the set of normal forms of elements in $H$ is a regular language. Part (1) of the theorem shows that all finitely generated subgroups are rational. This property directly shows that $G$ is Howson by an argument which is now well-known. If $H$ and $K$ have graphs $\Delta_2(H)$ and $\Delta_2(K)$ respectively, then the component of $\Delta_2(H) \times \Delta_2(K)$ which contains $(O_H, O_K)$ is the transition graph of a finite automaton which accepts a word $z$ in normal form if and only if $z \in H \cap K$, so $H \cap K$ is also rational. The Anisimov-Seifert theorem shows that any rational subgroup is finitely generated.

A subgroup $H$ is a *quasiconvex* if there exists a constant $d$ such that every point on a geodesic word which represents an element of $H$ is within distance $d$ of an element of $H$. Since a geodesic word is Dehn reduced, the graph $\Delta_2(H)$ is convex in the sense that any geodesic word representing an element of $H$ is the label on a path is $\Delta_2(H)$ which begins and ends at the basepoint $O_H$. Define the *diameter $d$* of $\Delta_2(H)$ to be the maximum distance of any vertex in the graph to the basepoint. It follows that $d$ is a quasiconvexity constant for $H$. Furthermore, the well-known Dijkstra Distance Algorithm runs in quadratic time and labels each vertex of $\Delta_2(H)$ by its distance from $O_H$, so $d$ is calculable in quadratic time.

The size of the graphs $\Delta_2(H)$ is "uniformly linearly bounded." Let

$$s_H = \sum_{i=1}^{m} |h_i|$$

be the sum of the lengths of the given generators of $H$. Let $k_G$ be the maximum length of any relator times the maximum of the $\rho_i$ (the number of relators involving $a_i$). The number of edges at any stage of the construction of $\Delta_2(H)$ will never exceed $k_G \, s_H$. We have implemented the calculation of $\Delta_2(H)$ with a computer program and for subgroups with $s_H$ around 1500 the calculation takes only a few seconds. Indeed, our ideas about the theorem have been formed to a great degree by computer experiments. The program will hopefully be available soon on the author's webpage.

## §2. The Construction.

To prove the theorem we now describe precisely the construction of the graph $\Delta_2(H)$ in two phases. Given generators $h_1, \ldots, h_m$ for $H$, we begin with the standard graph $\Delta_0(H)$ which is a bouquet of $m$ loops at a basepoint $O_H$ where the $i$-th loop is labelled by $h_i$. The number of edges in $\Delta_0(H)$ is exactly $s_H$.

A closed path in the subgroup graph which is labelled by a defining relator is called a *relator cycle.* If we see a path labelled by "enough" of a defining relator but which is *not* part of a relator cycle, we want to add the missing part of the relator to form a relator cycle. This enlarges the graph and we need the Reduction Hypothesis to guarantee that the process stops.

Here we use the beautiful idea of McCammond and Wise of perimeter reduction. We have chosen a subset $C \subseteq A$ such that every defining relator $(a_i a_j)^{m_{ij}}$ contains a generator from $C$. Given a current subgroup graph $\Delta$, for each edge $e$ labelled by a generator $a_i \in C$, we count the number of relator cycles which are missing at $e$. That is, the number of distinct relators $(a_i a_j)^{m_{ij}}$ such that $e$ is not on a relator cycle in $\Delta$ labelled by $(a_i a_j)^{m_{ij}}$.

The *count* $\gamma(\Delta)$ is the sum of four times the number of missing cycles over all edges labelled by generators from $C$ plus the current number of edges. For our purposes, "enough" of a relator is a relator with at most three letters missing. Our Reduction Hypothesis simply ensures that completing such paths to relator cycles reduces the count.

Suppose we see a path $\alpha$ in $\Delta$ labelled by a relator $(a_i a_j)^{m_{ij}}$ with three letters missing and $\alpha$ is not already on a relator cycle labelled by $(a_i a_j)^{m_{ij}}$. Then we add three new edges and two new vertices to the graph so that $\alpha e_1 e_2 e_3$ is a relator cycle.

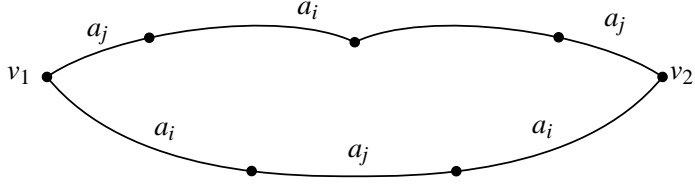

FIGURE 3

There are two cases in considering how this effects the count. First suppose that both $a_i$ and $a_j$ are in the chosen set $C$ of generators. The cycle $(a_i a_j)^{m_{ij}}$ is missing at each of the edges in $\alpha$, and we have now filled it in, so we have reduced the number of missing cycles by $2m_{ij} - 3$. The new edges are on the cycle we just formed so there are at most $3(\rho_{ij} - 1)$ missing cycles on the new edges. Since we multiply the number of missing cycles by 4 and we have added at most three new edges, the condition to reduce the count is just

$$2m_{ij} - 3 > 3(\rho_{ij} - 1)$$

which is $m_{ij} > \frac{3}{2}\rho_{ij}$.

The second case is that only one of the generators, say $a_i$, is in $C$ while $a_j \notin C$. If two of the new edges are labelled by $a_i$ then completing $\alpha$ to a relator cycle reduces the number of missing cycles by $m_{ij} - 2$ and increases the number of missing cycles by $2(\rho_i - 1)$. So the condition in this case is

$$m_{ij} - 2 > 2(\rho_i - 1)$$

or $m_{ij} > 2\rho_i$. If only one new edge is labelled by $a_i$ and the other two are labelled by $a_j$ the count is even further reduced. We alternate the following two operations.

(I) If there are two edges $e_1$ and $e_2$ with the same label $a_i$ incident at a vertex $v$, we fold the two edges. Note that folding reduces the number of missing cycles if $a_i \in C$ unless both $e_1$ and $e_2$ were already on all possible relator cycles involving the generator $a_i$. Even in this case or if $a_i \notin C$, folding reduces the number of edges so the count is reduced. We continue folding until the current graph is *trim,* that is, there are no pairs of edges to be folded.

(II) We search the current graph for a path $\alpha$ such that $\alpha$ is labelled by a relator with at most three letters missing but $\alpha$ is not part of a relator cycle. We then complete $\alpha$ to a relator cycle. Since $G$ satisfies the Reduction Hypothesis this operation reduces the count $\gamma$. Of course, the resulting graph need not be trim so we return to Step I.

We iterate (I) and (II) until neither is applicable. Since the original count for $\Delta_0(H)$ was at most $5k_G s_H$ and each step reduces the count and performing one step requires only searching once through the current graph, the process terminates with a graph $\Delta_1(H)$ in time quadratic in $s_H$.

For future reference, note that the iteration of Step II ensures that if there is a loop labelled by $(a_i a_j)^k$ at a vertex and $m_{ij} < \infty$ then this loop is part of a relator cycle so $k$ divides $m_{ij}$.

We again stress that we have used the Reduction Hypothesis to ensure that $\Delta_1(H)$ has the following property.

**The Relator Path Property.** *Any path $\alpha$ labelled by a relator with at most three letters missing is part of a relator cycle.*

Our construction of $\Delta_2(H)$ from $\Delta_1(H)$ depends only on this property and makes no further use of the Reduction Hypothesis. In complete generality, given any finite graph over $G$ which has the Relator Path Property, we can construct its 2-completion, so it is time to give a precise definition of what we mean. Since we also want to construct completions in the case of surface groups, we give a definition in the case of a group with a small cancellation presentation.

**2-Complete Graphs.** *Let $G = <X; R>$ be a finitely presented group where $R$ satisfies the small cancellation condition $C'(1/6)$ and all pieces have length one. A graph $\Gamma$ over $G$ is 2-complete if whenever two edges, $e_1 = (v_1, v)$ and $e_2 = (v, v_2)$ are incident at the vertex $v$ then $e_1 e_2$ is part of a relator cycle except possibly in the following situation: Both $v_1$ and $v_2$ have degree 2, and, there is an edge $e_3 = (v_3, v)$ incident at $v$ such that both $e_3 e_1^{-1}$ and $e_3 e_2^{-1}$ are part of relator cycles.*

*A 2-comletion $\Gamma_2$ of a graph $\Gamma$ is a finite 2-complete graph which contains $\Gamma$ as an embedded subgraph.*

Since we are stressing the independence of the two phases of our construction, we state the following as a theorem.

**Theorem II.** *Let $G$ be a Coxeter group of extra-large type. There is a uniform quadratic time algorithm which, when given any finite graph $\Delta_1$ over $G$ which has the Relator Path Property, constructs a 2-completion $\Delta_2$ of $\Delta_1$.*

*Proof.* We have to complete missing relator cycles. Call the edges of $\Delta_1$ *primary edges* and the vertices of $\Delta_1$ *primary vertices*. An $(i,j)$-*path* in $\Delta_1$ is a path of length at least two such that edge labels alternate between $a_i$ and $a_j$. The first edge may be labelled by either $a_i$ or $a_j$, so a $(j,i)$-path is the same thing as an $(i,j)$-path. Now $\Delta_1$ has the property that any maximal $(i,j)$-path which is not part of a relator cycle is a simple path and, since G is of extra-large type, requires adding at least four edges to complete it to a relator cycle.

For every maximal $(i,j)$-path in $\Delta_1$, $1 \leq i < j \leq n$, which is not already part of a relator cycle; we add the missing edges to complete that path to a relator cycle. New edges and vertices which are added in this process are called *secondary*. Call the resulting graph $\Delta_1'$.

Consider a primary vertex $v$ in $\Delta_1$ and the maximal paths beginning at $v$ which were completed to relator cycles by adding secondary edges. If such a path $\alpha$ is an $(i,j)$-path beginning with an edge labelled by $a_i$, then by the maximality of $\alpha$ there is no edge labelled by $a_j$ incident at $v$. There may also be an $(\ell, j)$-path beginning at $v$ whose first edge is labelled by $a_\ell$, $\ell \neq j, i$. See the figure below where the edges on the top of the figure are primary edges and the edges on the bottom are new secondary edges.

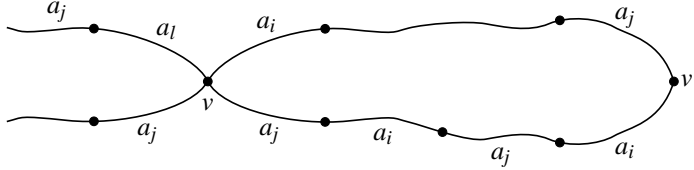

FIGURE 4

Identify all the secondary edges beginning at $v$ with label $a_j$, making the graph $\Delta_1'$ trim at the vertex $v$. Let $v'$ be the secondary endpoint of the secondary $j$-edge incident at $v$. We call such vertices *critical vertices*. Since the edges other than the $j$-edge at $v'$ are labelled by the distinct generators $a_i, a_\ell, \ldots$ beginning the maximal paths starting at $v$ which were completed to relator cycles, the new graph after folding the $j$-edges at $v$ is also trim at $v'$.

Hence, identifying new secondary edges with the same label which are incident at a primary vertex $v$ over all primary vertices $v \in \Delta_1$ results in a trim graph $\Delta_1''$ which is completely folded.

We still need to complete any pairs of edges, say labelled by $a_s$ and $a_t$ which are incident at a primary vertex $v$, and have $m_{st} < \infty$ but are not on an $(s,t)$-relator cycle. The construction above ensured that any two primary edges which meet at a primary vertex and are labelled by related generators lie on a relator cycle. But the current situation can occur when one or both of the edges of the pair are secondary edges.

First suppose that the edge $e_1$ labelled by $a_s$ is primary while the edge $e_2$ labelled by $a_t$ is secondary. Let $u$ be the other endpoint of $e_1$. There could also be a secondary edge $e_3$ labelled by $a_t$ incident at $u$. If this is the case complete $e_3 e_1 e_2$ to an $(s,t)$-relator cycle. If not, we complete $e_1 e_2$ to a relator cycle. If $e_1$ and $e_2$ are both secondary edges, complete $e_1 e_2$ to an $(s,t)$-relator cycle.

The result of having done this for all pairs of edges incident at primary vertices is $\Delta_2$. The claim is that $\Delta_2$ is trim as constructed. Take, for instance, the case where we complete two edges $e_1, e_2$ to an $(s,t)$-relator cycle as in the paragraph above. If there is a secondary edge $e_2$ labelled by $a_t$ incident at a primary vertex $v$ then there is no primary edge labelled by $a_t$ incident at that vertex. Also, $e_2$ was introduced in completing some maximal $(\ell, t)$-path to a relator cycle. So there is an $\ell$-edge incident at $v$. Since this cycle then exists, $s \neq \ell$. (There may be several choices of $\ell$ but none of them are equal to $s$.) So the new $s$-edge introduced at the end of $e_2$ keeps a trim graph. The argument is the same if there is another edge $e_3$ involved or in the case that both $e_1$ and $e_2$ are secondary.

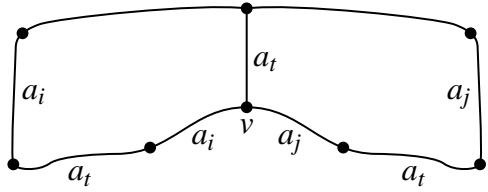

FIGURE 5

Recall that a *critical vertex* is a secondary vertex $v$ such that there is a secondary edge from $v$ to a primary vertex. The graph $\Delta_2$ has the crucial property that: *If edges labelled by $a_i$ and $a_j$ meet at any noncritical vertex then either $a_i$ and $a_j$ are not related or the edges lie on an $(i,j)$-relator cycle.* If edges $e_1, e_2$ labelled by $a_i$ and $a_j$ meet at a critical vertex $v$ and $a_i$ and $a_j$ are not on an $(i,j)$-relator cycle then there is a secondary edge labelled by some generator $a_t$, $t \neq i, j$ incident at $v$ and the other endpoints of $e_1$ and $e_2$ have degree 2 and the incident edges are also labelled by $a_t$.

To establish this last point we use the fact that since $G$ is of extra-large type all relator cycles newly introduced have at least four secondary edges. Only the secondary edges at each end of a path of secondary edges completing a relator cycle can also lie on other relator cycles. So if a secondary edge $e$ labelled by $a_i$ is incident to a critical vertex which is the endpoint of an edge labelled by $a_t$, the only edge incident at the other endpoint of $e$ is labelled by $a_t$. That all pairs of edges incident at a vertex and labelled by a related generators are on a relator cycle except in the one case discussed above is the crucial property of the 2-completion $\Delta_2$).

We have seen that $\Delta_2$ is trim. Any reduced path labelled by half a defining relator lies on a relator cycle. Indeed, the only path of length 2 not on a relator cycle is formed by edges labelled by $a_i$ and $a_j$ incident at a critical vertex as above and this is a maximal $(i,j)$-path. $\square$

Since 2-completion only involves completing relator cycles, it is clear that if we start with a subgroup graph $\Delta_1(H)$ with the Relator Path Property, its 2-completion $\Delta_2(H)$ contructed above is also a subgroup graph for $H$. The Normal Form Theorem thus shows that if $u$ is the label on a path from a vertex $v_1$ to a vertex $v_2$ in $\Delta_2(H)$ then there is a path from $v_1$ to $v_2$ labelled by the normal form $w$ of $u$. We say that $\Delta_2(H)$ is a *normal graph.* The stronger claim that $\Delta_2(H)$ is Dehn convex will be proved in the Appendix along with the Normal Form Theorem. As pointed out earlier, the diameter $d$ of $\Delta_2(H)$ can be calculated by Dijkstra's algorithm.

We now verify that $H$ has finite index if and only if $\Delta_2(H)$ is full, that is, there are edges labelled by all the generators of $G$ incident at all vertices. The "if" direction is clear. For each vertex $v$, pick a path $\alpha_v$ from the basepoint $O_H$ to $v$ and let $\hat{v}$ be the label on $\alpha_v$. If $w$ is any word on the generators of $G$, since $\Delta_2(H)$ is full, we read $w$ along a path starting at $O_H$ and ending at some vertex $v$. Thus $w(\hat{v})^{-1} \in H$ and $w$ is on the same right coset of $H$ as $\hat{v}$.

To show that if $\Delta_2(H)$ is not full then $H$ has infinite index, suppose that there is no edge labelled by a generator $a_\ell$ at some vertex $v$. Since $\Delta_2(H)$ is a normal graph, we can choose a path $\alpha_v$ from the basepoint $O_H$ to $v$ whose label is in normal form. We want to show that if $w$ is the label on the path $\alpha_v$ and $a_\ell$ is the generator missing at $v$ then $wa_\ell$ is still in normal form.

Suppose that $\alpha_v$ has length at least two. (The special cases where $\alpha_v$ is empty or has length one are easily dealt with on the same lines as the general argument.) Let $a_j$ be the label on the last edge $e_1$ in $\alpha_v$ and let $a_i$ be the label on the next to last edge $e_2$ in $\alpha_v$. Now $\ell \neq j$. If $v$ is primary then $i \neq \ell$ because $e_1 e_2$ lies on an $(i,j)$-relator cycle and $j$ is not missing at $v$. If $v$ is a critical secondary vertex there are two possible cases for the path $\alpha_v$ as pictured below.

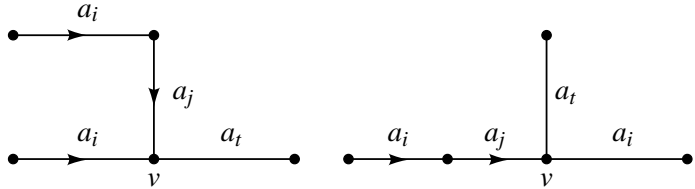

FIGURE 6

In both cases there must be an $a_i$ edge incident at $v$ so again $\ell \neq i$. If both endpoints of $e_1$ are secondary and not critical then $e_1$ lies only on an $(i,j)$-cycle and $\ell \neq i$. If the vertex on both $e_1$ and $e_2$ is critical and the last vertex of $e_1$ is a secondary noncritical vertex then it is possible that $\ell = i$ but in this case there is an edge labelled by some $a_t$, $t \neq i, j$ preceding $e_2$ in $\alpha_v$.

We can now construct a word $z$ such that for all $n \geq 1$, $z^n$ is in normal form and $z^n \notin H$. The existence of such an element shows that $H$ has infinite index. We have seen that $wa_\ell$ is in normal form and either $a_\ell$ is different from the last two generators in $w$ or $wa_\ell$ ends in $a_t a_i a_j a_i$ and $t$ is different from $i$ and $j$. Since there are at least three generators we can set $z = wa_\ell a_r a_s$ where $a_r$ differs from both $a_\ell$ and the last generator occurring $w$ and $a_s$ is different from $a_r$ and the first generator occurring in $w$. Then $z^n$ is in normal form for all $n \geq 1$ and since $\Delta_2(H)$ is a normal graph for $H$, $z^n \notin H$ for all $n \geq 1$.

The properties of $\Delta_2(H)$ are heavily used in the proof. Note, however, that if any relator cycles are added in Phase II then there are definitely secondary noncritical vertices of degree 2 and thus the graph is not full. So the graph is full if and only if Phase I perimeter reduction already results in a full graph.

How to decided membership in $H$ is now clear. Given an arbitrary word $w$, to decide if $w \in H$, calculate a Dehn reduced form $\hat{w}$ of $w$ and then read $\hat{w}$ on the graph $\Delta_2(H)$ and see if one stays in the graph and returns to the basepoint at the end of $\hat{w}$.

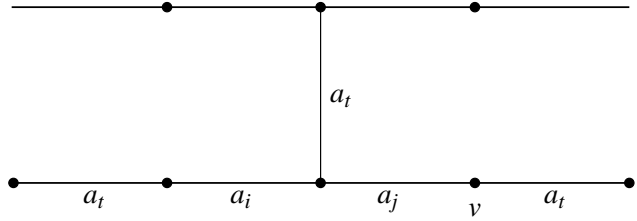

Figure 7

We now turn to separability. If $w \notin H$, we modify $\Delta_2(H)$ to $\Delta_0(H, w)$ by adding a stem to $\Delta_2(H)$ which is a path consisting of new edges labelled by $w$ starting at the basepoint $O_H$ and having terminal vertex $T_w$. Then we again carry out both perimeter reduction and 2-completion to obtain a graph $\Delta_2(H, w)$ having the same properties as $\Delta_2(H)$ but with a second distinguished vertex $T_w \neq O_H$ and a path from $O_H$ to $T_w$ with label $w$.

We now show that if $G$ also satisfies the Separability Condition we can directly define an action of $G$ on the set $V$ of vertices of $\Delta_2(H, w)$ such that $O_H = O_H \cdot g$ if and only if $g \in H$ but $T_w = O_H \cdot w$. The corresponding homomorphism $\varphi$ thus has $\varphi(w) \notin \varphi(H)$. Since all generators of $G$ have order 2, there is a natural candidate for such an action. Namely, if there is an edge labelled by $a_i$ in $\Delta_2(H, w)$ incident at a vertex $v$ and the other endpoint of the edge is $v_1$ (possibly $v_2 = v$) then $a_i$ interchanges $v$ and $v_1$. If there is no edge incident at $v$ with label $a_i$ then $a_i$ fixes $v$. This definition has the property that $v = v \cdot a_i^2$ for all $i$ and all $v$.

We now have to check that if $a_i$ and $a_j$ are any two distinct generators which are related and $v$ is any vertex then $v = v \cdot (a_i a_j)^{m_{ij}}$ and we really have an action of the Coxeter group $G$. If there is neither an $i$-edge nor a $j$-edge incident at $v$ then $a_i$ and $a_j$ both fix $v$ so $v \cdot (a_i a_j)^k = v$ for all $k$. Suppose first that $v$ is a primary vertex. If there are both an $i$-edge and a $j$-edge incident at $v$, then there is an $(i, j)$-relator cycle beginning at $v$ and thus $v \cdot (a_i a_j)^{m_{ij}} = v$. If there is, say, an $i$-edge $e_1$ incident at $v$ but not a $j$-edge, then there is also not a $j$-edge incident at the other endpoint $v_1$ of $e_1$ by the properties of $\Delta_2(H, w)$. (If there were then $e_1$ would again be on an $(i, j)$-relator cycle.)

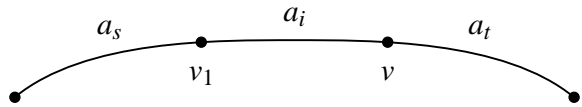

Figure 8

Now $v \cdot (a_i a_j)^2 = v$ since

$$v \cdot a_i a_j a_i a_j = v_1 \cdot a_j a_i a_j = v_1 \cdot a_i a_j = v \cdot a_j = v.$$

Thus the assumption that $m_{ij}$ is even ensures that $v \cdot (a_i a_j)^{m_{ij}} = v$.

Similarly, if $v$ is a noncritical secondary vertex, the only two edges incident at $v$ lie on a relator cycle. If this is an $(i, j)$-relator cycle, $v \cdot (a_i a_j)^{m_{ij}} = v$. If there is only an $i$-edge $e_1$ incident at $v$ and the other endpoint $v_1$ of $e_1$ is not a critical vertex, then there is again no $j$-edge incident at $v_1$ and the argument of the preceding paragraph applies.

$\Delta_2(H, w)$ has the property that the only kind of vertex at which both an $i$-edge and a $j$-edge can meet and not lie on a relator cycle is a critical secondary vertex $v$.

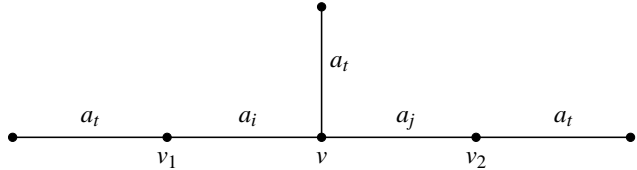

FIGURE 9

In this case there is also a $t$-edge, $t \neq i, j$ incident at $v$ and $a_i$ and $a_j$ are both related to $a_t$. (Thus $a_i, a_j$ and $a_t$ lie on a triangle in the modified Coxeter graph $\Gamma$.) There is no $j$-edge incident at $v_1$ and there is no $i$-edge incident at $v_2$. Now we have $v \cdot (a_i a_j)^3 = v$. Since $a_i$ and $a_j$ lie on a common triangle $v$ in $\Gamma$, the Separability Condition ensures that $m_{ij}$ is divisible by 3, so $v \cdot (a_i a_j)^{m_{ij}} = v$. The reader can easily verify that in this situation we also have $v_1 \cdot (a_i a_j)^3 = v_1$ and $v_2 \cdot (a_i a_j)^3 = v_2$, which takes care of the last remaining case.

As pointed out in the introduction, many Coxeter groups are not subgroup separable since there are Coxeter groups with unsolvable membership problems. Proving subgroup separability for Coxeter groups which have small exponents seems to require a great deal of specific geometric information. See for example [L-R].

Perimeter reduction gives a nice proof that the Reduction Hypothesis ensures that the Phase I construction of $\Delta_1(H)$ terminates. Extensive computer investigation very strongly suggests, however, that Phase I always terminates for any finitely generated subgroup of a Coxeter group of extra-large type but proving this would seem to be combinatorially complicated. We believe the following conjecture.

**Conjecture.** *Theorem I is true as stated for all Coxeter groups of extra-large type.*

We certainly do not have a good feeling for how essential the Separability Condition is with regards to actual subgroup separability but it is perhaps illuminating to consider residual finiteness. It is well-known that every finitely generated Coxeter group $G$ is residually finite since $G$ is a finitely generated subgroup of some $GL(n, \mathbb{R})$. In general however, given a presentation and an element $w$ the linear group argument does not at all tell us how to actually find a homomorphism from $G$ to a finite group with $\varphi(w) \neq 1$. If $G$ is of extra-large type and satisfies the Separabilty Condition we have the following.

**Theorem III.** *There is a uniform quadratic time algorithm which, when given any Coxeter group $G$ of extra-large type satisfying the Separability Condition and any element $w \neq 1$ of $G$. writes out a homomorphism $\phi$ from $G$ to a finite symmetric group $S$ such that $\phi(w) \neq 1$.*

*Proof* Given a nontrivial element of $G$, let $w$ be a Dehn reduced form. Let $\Delta_1(w)$ be the graph consisting of a path of successive edges labelled by the letters in $w$. Since $w \neq 1$, the initial vertex $O_w$ and the terminal vertex $T_w$ of $\Delta_1(w)$ are distinct. Since $w$ is Dehn reduced it does not contain more than half of a defining relator and thus $\Delta_1(w)$ has the Relator Path Property that there are no $(i,j)$-paths requiring fewer than four edges to be completed to a relator cycle. We can then construct $\Delta_2(w)$ from $\Delta_1(w)$ as before. Since $G$ satisfies the Separability Condition , we can define the action of $G$ on the set $v$ of vertices of $\Delta_2(w)$ as before and $w$ sends $O_w$ to $T_w$ so $w$ is not in the kernel of the action. □

Recall that Peter Scott [SC] proved that all Fuchsian groups are subgroup separable by showing that surface groups are subgroup separable and these are separable because an appropriate Coxeter group is subgroup separable. Scott showed that subgroup separability is a *virtual property*: If $[K : G] < \infty$ and $K$ is subgroup separable then $G$ is subgroup separable. So any finitely generated group containing a subgroup separable Coxeter group of finite index is itself subgroup separable.

## §3. Surface Groups.

McCammond and Wise pointed out that surface groups form an important class of groups to which we can apply perimeter reduction. We point out here that constructing the 2-completion of a subgraph graph works for surface groups in exactly the same way as for Coxeter groups. The 2-completion therefore gives all the detailed information of our Theorem except for subgroup separability. The Reduction Hypothesis we need is that the standard defining relator has length at least 8, that is, the genus is at least 2 in the orientable case and is at least 4 in the nonorientable case. As usual, the symmetrized set $R$ of defining relators consists of all cyclic permutations of the standard relator and its inverse.

We briefly discuss the construction. Given a finite set $\{h_1, \ldots, h_m\}$ of generators of a subgroup $H$ of a surface group $G$ satisfying the Reduction Hypothesis, we begin with the bouquet of circles $\Delta_0(H)$ as always. Now, of course, the generators have infinite order.

Phase I now consists of completing paths which are labelled by more than half of a defining relator to relator cycles. Let $h$ be half the length of the defining relator. The *count* $\gamma(\Delta)$ for a graph $\Delta$ is the sum over all edges $e$ of $\Delta$ of $h$ times the number of cycles missing at $e$ plus the total number of edges.

Suppose that we complete a path $\alpha$ of length $\ell > h$ to a relator cycle. We add $k < h$ new edges. Because the defining relator is quadratic each generator begins only two distinct relator cycles. Thus we have removed $\ell$ missing cycles and added at most $k$ missing cycles and the count is reduced. So completing paths containing more than half of a defining relator and folding terminates in a linear number of steps by perimeter reduction, yielding the graph $\Delta_1(H)$. As in the case of Coxeter groups, reading the two initial letters of an element $r$ of the symmetrized set $R$ uniquely determines $r$ and any path not already part of a relator cycle requires at least four edges to complete. The construction of the 2-completion $\Delta_2(H)$ thus proceeds exactly as in the case of Coxeter groups.

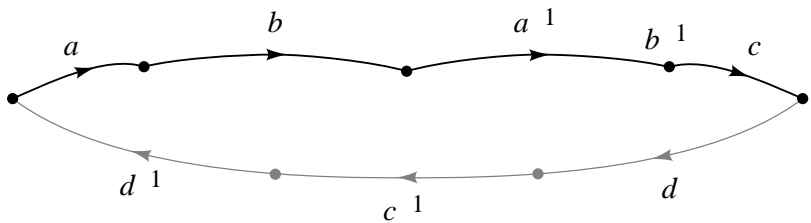

Figure 10

As before, $\Delta_2(H)$ is the transition graph of a finite automaton which accepts a Dehn reduced word $w$ if and only if $w \in H$ and the diameter of $\Delta_2(H)$ is a quasiconvexity constant for $H$. Also as before, $H$ is of finite index if and only if $\Delta_2(H)$ is a full graph.

Pittet[P] earlier showed that fundamental groups of orientable surfaces of genus at least two are locally quasiconvex. The argument of McCammond and Wise [M-W] and our argument here show that such surface groups are locally quasiconvex and have uniform membership problem solvable in quadratic time. In his thesis [G] Grunschlag proved the following theorem.

**Theorem** (Grunschlag). *If $G$ is a finitely generated group whose uniform membership problem is solvable in time $f(n)$, and $J$ contains $G$ as a subgroup a subgroup of finite index, then the uniform membership problem for $J$ is solvable in time $f^2(n)$.*

In particular, this theorem shows that having uniform membership problem solvable in polynomial time is a virtual property. Also, several authors, see for example Kapovich and Short[KS2], have shown that local quasiconvexity is a virtual property.

**Theorem IV.** *If $F$ is any group containing as a subgroup of finite index either an orientable surface group of genus at least two or a Coxeter group which satisfie the Reduction Hypothesis then $F$ is locally quasiconvex and the uniform membership problem for $F$ is solvable in at most quartic time.*

In particular, the theorem applies to finitely generated discrete groups of isometries of the hyperbolic plane. It seems to the author that having uniform membership problem solvable in polynomial time is an interesting property which may help point out which groups have "really good" geometry.

## §4. Appendix
## Normal Forms and Dehn Convexity

In this appendix we sketch a proof of the facts we used about normal forms and geodesics. Let $G = \langle X; R \rangle$ be a finitely presented group where $R$ is a symmetrized set of relators satisfying the small cancellation condition $C'\left(\frac{1}{7}\right)$ (See [L-S] for the basic definitions of small cancellation theory.) For convenience we suppose that all pieces have length 1 since this is the case in the groups which we have considered. We want to prove that the subgroup graphs $\Delta_2(H)$ which we have constructed are indeed Dehn convex in the sense that if there is a path $\alpha$ between two vertices labelled by a reduced word $w$ and $z$ is a Dehn reduced word with $z = w$ in $G$ then there is a path labelled by $z$ between the two vertices.

Our $\Delta_2(H)$ all have the property that any path $\gamma$ labelled by a subword of an element of $R$ containing at least three letters is part of a relator cycle. So let $\Gamma$ be a finite graph over the group $G$ which has this property. We shall show that $\Gamma$ is Dehn convex. Let $w$ be the label on a reduced path between two fixed vertices. Certainly, if $w$ contains a subword $s$ which is more than half a defining relator, the part of $w$ labelled by $s$ lies in a relator cycle and we can perform a Dehn reduction by going around the shorter part of the relator cycle and obtain a shorter path between the two vertices. Iterating this procedure, we may assume that $w$ is $R$-reduced.

Let $z$ be the shortlex normal form of $w$ to show that there is a path in $\Gamma$ labelled by $z$ between the two vertices. We consider the *equality diagram* $E(w, z)$ of $w$ and $z$, which is formed as follows. Write $w \equiv pw'$ and $z \equiv pz'$ in the free group where $p$ is the maximal common initial segment of $w$ and $z$. The prefix $p$ may be empty but if either $w'$ or $z'$ is empty then $w \equiv z$ is the free group. (If we had $z'$ empty and $|w'| > 0$ we would have a nontrivial $R$-reduced word $w'$ equal to 1 in $G$ which is impossible.) Now write $w \equiv pw_1 d$, $z \equiv pz_1 d$ where $d$ is the largest common terminal segment of $w'$ and $z'$. Again, neither $w_1$ nor $z_1$ are empty unless $w \equiv z$ and is thus already in normal form and we are done.

Consider the minimal reduced diagram $C$ for the nontrivial cyclically reduced word $w_1 z_1^{-1}$. The equality diagram $E(w, z)$ is the reduced $R$-diagram which reads $w$ along the top boundary $\sigma$ and $z$ along the bottom boundary $\tau$. (So a boundary cycle of $E(w, z)$ is $\sigma\tau^{-1}$.) Equality diagrams have been often used in small cancellation theory and have a very simple structure which seems to first appear in Schupp's [S1] analysis of the structure of conjugary diagrams for certain small cancellation groups. Since no region can have more than half of its boundary as a consecutive part of either $\sigma$ or $\tau$, successive applications of Greendlinger's Lemma show that each region $D$ has interior degree at most 2 and the intersection of $\partial D$ with the top (bottom) boundary of $E(w,z)$ is connected. Thus $E(w,z)$ of some number of "islands" connected by segments or meeting at a single vertex.

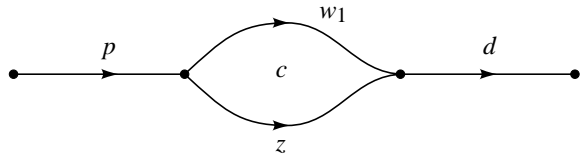

FIGURE 11

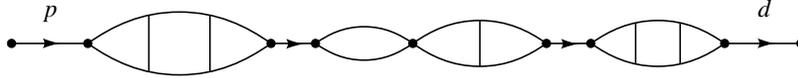

FIGURE 12

In each island, a region has interior degree 0 if it is the entire island, interior degree 1 if it is at the end of an island and interior degree 2 otherwise. Since both $w$ and $z$ are $R$-reduced and $R$ satisfies $C'\left(\frac{1}{7}\right)$, for each region $D$, the segment which comprises the intersection of $\partial D \cap \sigma$ (respectively $\partial D \cap \tau$) must contain at least two letters and is thus labelled by a nonpiece.

A vertex of $\sigma$ and a vertex of $\tau$ which are connected by an interior edge are said to *correspond*. A vertex in $\sigma \cap \tau$ corresponds to itself. Now the word $w$ is both the label on the path $\alpha$ in $\Gamma$ between vertices $v_0$ and $v_n$ and the label on the top boundary $\sigma$ of the equality diagram $E(w,z)$. Working from left to right in $E(w,z)$ we can successively show that for each initial segment $u$ of $w$ which ends at a vertex $s_i$ corresponding to a vertex $t_i \in \tau$, there is a path $\beta_i$ in $\Gamma$ between $v_0$ and the vertex $v_i$ at which $\beta_i$ ends in $\Gamma$ which has the same label as the path in $E(w,z)$ obtained by starting at its initial vertex $s_0$, reading along $\tau$ until the vertex $t_i$ corresponding to $s_i$ is encountered and then going along the interior edge from $t_i$ to $s_i$. (The latter provided $t_i \neq s_i$.) The picture is

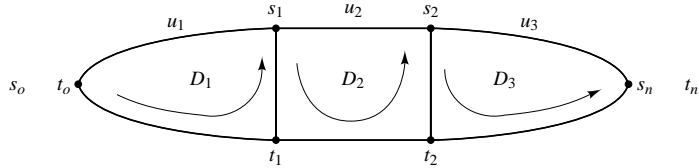

FIGURE 13

Suppose, for example, that $E(w,z)$ begins exactly as in the figure above. The word $u_1$ must contain at least three letters of the relator cycle labelling the region $D_1$. Thus the subpath of $\alpha$ labelled by $u_1$ lies on a relator cycle in the graph $\Gamma$ and going around that cycle the other way gives the desired path $\beta_1$ in $\Gamma$ from $v_0$ to vertex $v_i$. The only point to notice is that in any event the label on a segment $\partial D_i \cap \sigma$ contains at least two letters. Thus the path in $\Gamma$ starting with the last edge in $\beta_1$ and then following that part of $\alpha$ labelled by $s_2$ must contain at least three letters. This path thus is also part of a relator cycle and we can go around the other way. Following $\beta_1$ by this path and omitting the edge traversed twice we obtain a reduced path $\beta_2$ in $\Gamma$ from $v_0$ to the vertex $1/2$ labelled by the desired word obtained by reading along $\tau$ from $O_w$ to $t_2$ and then taking the interior edge to $s_2$. We can clearly continue in this way, finally obtaining the desired reduced path $\beta$ with label $z$ which goes from $v_o$ to $v_n$ in $\Gamma$.

We have thus shown that if $w$ is the label on a path in $\Gamma$ between two vertices then there is a path between these vertices whose label is the normal form $z$ of $w$. Now let $w'$ be any $R$-reduced word equal to $z$ in the group $G$. We now take the equality diagram $E(w', z)$ and work the other way. We have the path labelled by $z$ existing in $\Gamma$. Apply the argument above but work from $\tau$ to $\sigma$ to construct a path labelled by $w'$ in $\Gamma$. In particular, this shows that $\Gamma$ is Dehn convex in the sense that if $w$ is a word labelling a path $\alpha$ in $\Gamma$ and $u$ is any Dehn reduced word equal to $w$ in $G$ then there is a path labelled by $u$ between the endpoints of $\alpha$.

The argument of working from left to right in an equality diagram also shows that the set of normal forms of elements of $G$ is a regular language. First of all, the set of all Dehn reduced words is certainly regular since it is defined by not containing a subword from a fixed finite list. Our automaton verifies that an input $w$ is Dehn reduced and as long as it keeps track only of one possible region in an equality diagram plus two pieces of information. For example, suppose that $G$ is the Coxeter group

$$G = \langle a_1, a_2, a_3;\ a_i^2,\ (a_i a_j)^4,\ i \neq j \rangle$$

The point is that at any point in a Dehn reduced word $w$ the possible corresponding part of the other boundary $\tau$ of $E(w, z)$, $z$ a normal form, is completely determined. Suppose, for example, that $w = a_2 a_1 a_2 a_1 a_3 a_2 a_3 a_1 a_2 a_1 a_2 a_3$. When two consecutive letters do not form part of a region, we start a new region and we start a new region at the beginning of $w$. Once two letters are part of a relator, that relator is uniquely determined and thus the label in the bottom part. On reading $a_2 a_1$ we know that a first region starts

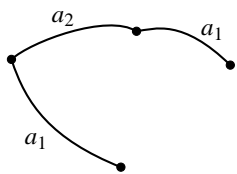

Figure 14

When we arrive at a break in $w$ of the current relator, we know how to complete the first region. In our example we have

The two pieces of information which the automaton needs to remember is whether the proposed label for the normal form is lexicographically less than the word being read, in this example it is, and the difference in length between the top and the bottom. Notice that this length cannot be greater than 1. If it is, we can find a shorter word equal to $w$ by going around the bottom and then traversing the interior

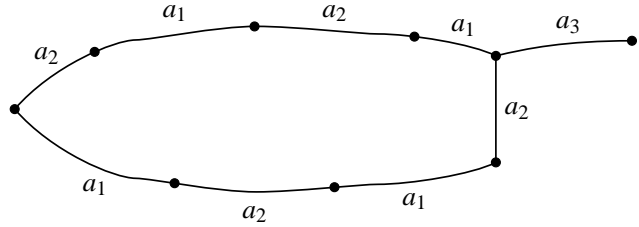

FIGURE 15

edge up to the corresponding vertex of $w$. In the example, when we read $a_3$, we switch to a new region and at the end of reading the segment $a_3 a_2 a_3$ we have the picture

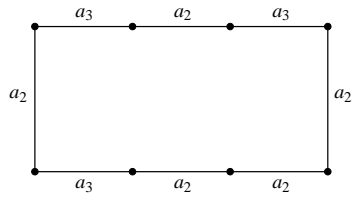

FIGURE 16

The bottom remains 1 letter shorter than the top. On reading the last letter $a_3$, we have the picture

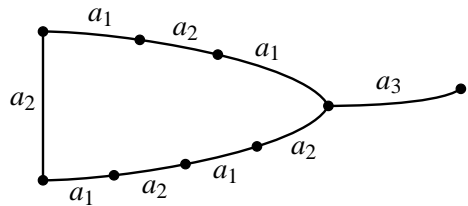

FIGURE 17

The top and bottom now haved the same length but the bottom label comes first lexicographically so $w$ is rejected at this point. Since there are only finitely many "region pictures" we can use them as the states of a finite automaton and we are done. Technically we need a special end marker symbol to indicate that we have come to the end of the word, but if the language with the special endmarker is regular then so is the language without the endmarker by standard closure properties of regular languages.

Paul E. Schupp
Department of Mathematics
University of Illinois
Urbana, IL 61801
schupp@math.uiuc.edu

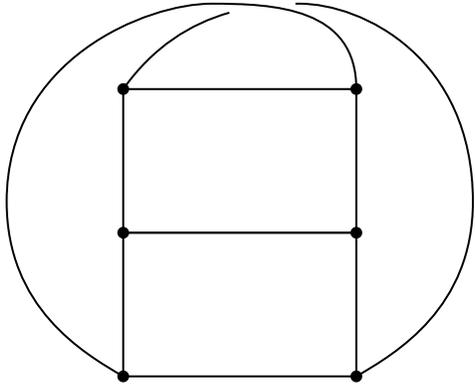